\documentclass[12pt,leqno]{amsart}

\usepackage{enumerate}
\usepackage{float} 
\usepackage{hyperref}
\usepackage{graphicx,subfigure}
\usepackage{amsmath,amsthm,amssymb,amsfonts}%,hypernat}
\theoremstyle{plain}
\newtheorem{theorem}{Theorem}
\newtheorem{lemma}[theorem]{Lemma}
\newtheorem{proposition}[theorem]{Proposition}
\newtheorem{corollary}[theorem]{Corollary}

\theoremstyle{definition} 
\newtheorem{remark}{Remark}
\floatstyle{boxed} \newfloat{alg}{htb}{}
\floatname{alg}{Algorithm}
\newcommand{\R}{{\mathbb R}}
\newcommand{\N}{{\mathbb N}}

\renewcommand{\rho}{{\varrho}}

\def\e{{\varepsilon}}

\def\phi{{\varphi}}  % huebscher!? 

\renewcommand{\rho}{{\varrho}}

\newcommand{\abs}[1]{\left\vert #1 \right\vert} 
\newcommand{\norm}[2]{\left\Vert #1  \right\Vert _{#2}} 
 
\newcommand{\set}[1]{\left\{#1\right\}}
\newcommand{\expect}{\mathbf E}

\newcommand{\scalar}[2]{\left\langle #1,#2\right\rangle}

\begin{document}

\title[Error bounds of MCMC]
{Error bounds for computing the expectation by Markov chain Monte Carlo} 

\author{Daniel Rudolf} 
\thanks{This work was supported by the DFG Priority Program 1324.}
\address{ Friedrich Schiller University Jena, 
Mathem. Institute,
Ernst-Abbe-Platz 2, 
D-07743 Jena, Germany}
\email{daniel.rudolf@uni-jena.de}
\date{Version: \today}

\begin{abstract}
We study the error of reversible Markov chain Monte Carlo methods for approximating the expectation of a function.
Explicit error bounds with respect to the $\ell_2$-, $\ell_4$- and $\ell_\infty$-norm of the function are proven.
By the estimation the well known asymptotical limit of the error is attained, i.e. our bounds are correct to first order as $n\to\infty$.
We discuss the dependence of the error on a burn-in of the Markov chain.
Furthermore we suggest and justify a specific burn-in for optimizing the algorithm. 
\end{abstract}

\maketitle

% \begin{document}

%*******************************%
% Please input your paper here. %

\section{Introduction}

We start with
a probability distribution $\pi$ on a finite set $D$ and a function $f: D\to \R$. 
The goal is to compute the expectation denoted by
\[
S(f)=\sum_{x\in D} f(x)\pi(x).
\]
Let the cardinality of $D$ be very large such that an exact computation of the sum is practically impossible.
Furthermore suppose that the desired distribution is not explicitly given, i.e. we have no random number generator for $\pi$ available.
Such kind of problems arise in statistical physics, in statistics, and in financial mathematics (see for instance \cite{MC_MC_pract,liu}).
The idea of approximating $S(f)$ via Markov chain Monte Carlo (MCMC) is the following:
Run a Markov chain on $D$ to simulate the distribution $\pi$ and compute the time average over the last $n$ steps. % of the chain.
Let $X_1,\dots,X_{n+n_0}$ be the chain, %denotes the sample
then we obtain as approximation
\[
S_{n,n_0}(f)=\frac{1}{n}\sum_{i=1}^n f(X_{i+n_0 }).
\] 
By $n_0$ the so called burn-in is given, loosely spoken this is the number of time steps 
taken to warm up. Afterwards the distribution of the generated Markov chain is 
(hopefully) close to the stationary one.\\
%taken to get close to the stationary distribution.
%such that the distribution of the generated Markov chain is close to the desired one. 

A Markov chain is identified with its initial distribution $\nu$ and its transition matrix $P$. We
restrict ourself to reversible chains which are ergodic, i.e. the second largest absolute value $\beta$ of the eigenvalues of $P$ is 
smaller than one.
%It is well known that one may use this to consider how fast the distribution of the chain reaches stationarity 
It is well known that the distribution of these chains reaches stationarity exponentially %(in $O(\beta^n)$ with $n\to \infty$) 
(see \cite{bremaud,rosenthal_hybrid,peres}). \\

The error of $S_{n,n_0}$ 
for $f\in\R^D$ is measured by
\[
e_\nu(S_{n,n_0},f)=\left(\expect_{\nu,P}\abs{S_{n,n_0}(f)-S(f)}^2\right)^{1/2},
\]
where $\expect_{\nu,P}$ denotes the expectation of the Markov chain. % which is determined by $P$ and initial distribution $\nu$.\\
The asymptotic behavior of the integration error can be written in terms of the eigenvalues and eigenfunctions of $P$.  It holds true that
\[
\lim_{n\to\infty}n\cdot e_\nu(S_{n,n_0},f)^2 \leq \frac{1+\beta_1}{1-\beta_1}\norm{f}{2}^2, 
\]
where $\beta_1$ is the second largest eigenvalue (see \cite{sokal_lect,mathe1}).
The constant $\frac{1+\beta_1}{1-\beta_1}$ is optimal but this statement does not give an error bound for finite $n$ and also does not include anything concerning the choice of $n_0$.
How does an explicit error bound of the MCMC method look like 
where the asymptotic behavior is attained?\\

Let us give an outline of the structure and the main results. 
Section~\ref{prel} contains the used notation and presents some relevant statements concerning Markov chains.  Section~\ref{err_sec} contains the new results.  %We link this results to convergence properties of the Markov chain to get upper error bounds.  
The explicit error bound is developed with respect to the $\ell_2$-, $\ell_4$- and $\ell_\infty$-norm of the function $f$.
%Let us illustrate exemplary the bound with adapted burn-in for $\norm{f}{\infty}\leq1$. 
For $\norm{f}{\infty}\leq1$ and $C=2 \sqrt{\norm{\frac{\nu}{\pi}-1}{\infty}}$ we obtain the following. The error obeys
\[
e_\nu(S_{n,n_0},f)^2 \leq \frac{2 }{n(1-\beta_1)}+\frac{2C\beta^{n_0}}{n^2(1-\beta)^2}.
\]
For details and estimates concerning $\ell_2$ and $\ell_4$ we refer to Theorem~\ref{main_thm} in Section~\ref{main_thm_sec}. In Section~\ref{burn_in_sec} it turns out that $n_0 = \max\set{\left\lceil \frac{\log\left(C\right)}{\log(\beta^{-1})} \right\rceil,0}$ is a reasonable choice for the burn-in. 
Then the error bound simplifies to
\[
e_\nu(S_{n,n_0},f)^2 \leq \frac{2 }{n(1-\beta_1)}+\frac{2}{n^2(1-\beta)^2}.
\]
%%%%%%%%%%%%
%%%%%%%%%%%%
%%%%%%%%%%%%
For an interpretation let us consider the case where $\beta=\beta_1$. 
Then the cost $n+n_0$ which are needed for an optimal algorithm 
to approximate $S(f)$ within an error of $\e$ can be bounded by
\[
	 \left\lceil \frac{4}{\e^2 (1-\beta)} \right\rceil 
											 + \left\lceil \frac{\log(C)}{\log(\beta^{-1})} \right\rceil.
\]
Hence $C$ can be astronomically large, for instance depending polynomially on the cardinality of the state space,
because it comes in logarithmically.\\

%%%%%%%%%%%%
%%%%%%%%%%%%
%%%%%%%%%%%%
In many examples a good estimate for $\beta$ can be achieved, see for instance \cite{randall_decomp, bassetti_dia, bass_leisen}.
Therefore it is straightforward to apply the explicit error bound.

\section{Preliminaries}%Definitions and Notations} 
\label{prel}

The Markov chain $X_1,X_2,\dots$ is a stochastic process with state space $D$. 
It is given by initial distribution $\nu$ %i.e. %$\nu=(\nu(1),\dots,\nu(m))$ and 
%$\sum_{x\in D}\nu(x)=1$
and
transition matrix $P=(p(x,y))_{x,y \in D}$ and denoted by $(\nu,P)$.
For $x,y\in D$ the entry $p(x,y)$ presents the probability of jumping from state $x$ 
to state $y$ in one step of the chain.\\

 %In detail, it is the expectation with respect to the probability
%\[
%W_{\nu,P}(x_1,\dots,x_{n+n_0})=\nu(x_1) p(x_1,x_2)\cdot\dots\cdot p(x_{n+n_0-1},x_{n+n_0}),
%\]
%which is defined on $(D^{n+n_0},\mathcal{P}(D^{n+n_0}))$.
%\noindent
%If we apply $P$ from the left to a function $f$ 
By $Pf(x)=\sum_{y\in D} p(x,y) f(y)$ we obtain the expectation of the value of $f\in\R^{D}$ after one step of the chain starting from $x\in D$. The expectation after $k$ steps of the Markov chain from $x$ is given by 
$P^kf(x)=\sum_{y\in D} p^k(x,y) f(y)$, where $P^k=(p^k(x,y))_{x,y \in D}$ denotes the $k$-th power of $P$. 
Similarly we consider the application of $P$ to a distribution $\nu$, i.e. $\nu P(x)= \sum_{y\in D} \nu(y) p(y,x)$. This is the distribution which arises after one step where the initial state was chosen by $\nu$. The distribution which arises after $k$ steps is given by $\nu P^k(x)= \sum_{y\in D} \nu(y) p^k(y,x)$.
\\

The expectation $\expect_{\nu,P}$ of the Markov chain $X_1,\dots,X_{n+n_0}$ is taken with respect to the probability measure
\[
W_{\nu,P}(x_1,\dots,x_{n+n_0})=\nu(x_1) p(x_1,x_2)\cdot \dots \cdot p(x_{n+n_0-1},x_{n+n_0}),\quad n,n_0\in \N,
\] 
on $D^{n+n_0}$.
Using this for $i\leq j$ we obtain a characterization by the transition matrix 
\begin{equation}
  \label{err_help}
	\expect_{\nu,P}({f(X_i)f(X_j)})=\sum_{x\in D} P^i(f P^{j-i}f)(x) \nu(x).
\end{equation}

\subsection{Reversibility and spectral structure}
We call the Markov chain with transition matrix $P$, or simply $P$, reversible with respect to a probability measure $\pi$ if the detailed balance condition
\[
\pi(x) p(x,y)=\pi(y) p(y,x)
\] 
holds true for $x,y \in D$. 
If $P$ is reversible, then $\pi$ is called stationary distribution of the Markov chain, i.e. $\pi P(x)=\pi(x)$.
%Since $P$ is reversible, $\pi$ is a stationary distribution of the Markov chain, i.e. $P\pi(x)=\pi(x)$. 
Note that, if $P$ is reversible then $P^k$ is also reversible.
%Furthermore we always assume that the Markov chain is irreducible and aperiodic (for an introduction see \cite{bremaud}) therefore the stationary distribution $\pi$ of $P$ is unique. 
%The next outcome describes a property of reversible chains and follows for discrete state spaces evidently from the definition of reversibility. For continuous state space see \cite{expl_error}. 
%
%\begin{lemma}
%Let $P$ be a reversible transition matrix and let $F:D\times D\to\R$. Then we have for
%$P^k=(p^k(x,y))_{x,y \in D}$ denoting the $k$-th power of $P$ that 
%\begin{equation}
%	\sum_{x\in D} \sum_{y\in D} F(x,y) p^k(x,y) \pi(x) = \sum_{x\in D} \sum_{y\in D} F(y,x) p^k(x,y) \pi(x).
%	\label{F_s_l_i}
%\end{equation}
%\end{lemma}
Let us define the weighted scalar-product
\[
\scalar{f}{g}_\pi=\sum_{x\in D} f(x)g(x)\pi(x),
\]
for functions $f,g\in\R^D$. Then let $\norm{f}{2}=\scalar{f}{f}_\pi^{1/2}$.
By considering the scalar-product it is easy to show, that reversibility is equivalent to $P$ being self-adjoint. 
%Let $\norm{f}{2}=\scalar{f}{f}_\pi^{1/2}$ be the (weighted) norm of $f\in\R^D$. 
Furthermore suppose that the underlying Markov chain is irreducible and aperiodic, 
this is also called ergodic. 
For details of these conditions we refer to the literature, for instance \cite{haggstrom,bremaud,peres}. 
It is a well known fact that this implies the uniqueness of the stationary distribution. % $\pi$ of $P$. 
Applying the spectral theorem of self-adjoint stochastic matrices and ergodicity we obtain that $P$ has real eigenvalues 
\[
1=\beta_0>\beta_1\geq \beta_2 \geq \dots \geq \beta_{\abs{D}-1} > -1 
\]
with a basis of orthogonal eigenfunctions $u_i$ for $i\in\set{0,\dots,\abs{D}-1}$, i.e.
\[
Pu_i=\beta_i u_i, \qquad \scalar{u_i}{u_j}_\pi = \delta_{ij} = \begin{cases} 1 & i=j\\
																												            0 & i\not=j.
																												            \end{cases}
\]
%Note that, $ \beta:=\max \set{\beta_1,\abs{\beta_{\abs{D}-1}}}<1 $ is an implication of the irreducibility and aperiodicity.
Additionally one can see that $u_0(x)=\mathbf{1}$ and $S(u_i)=0$ for $i>0$. 

\subsection{Convergence of the chain}
\label{convergence}
The speed of convergence of the Markov chain to stationarity is measured
%The total variation norm provides one tool for doing so i.e.
%let be $\nu$, $\mu$ distributions on $D$ then 
%\[
%\norm{\nu(\cdot)-\mu(\cdot)}{\text{tv}}:=\frac{1}{2}\sum_{x\in D} \abs{\nu(x)-\mu(x)}.
%\] 
%Note, that always $\norm{\nu(\cdot)-\mu(\cdot)}{\text{tv}}\leq1$.
by the so called $\chi^2$-contrast. Let $\nu$, $\mu$ be distributions on $D$ then 
\[
\chi^2(\nu,\mu)=\sum_{x \in D} \frac{(\nu(x)-\mu(x))^2}{\mu(x)}.
\]
The $\chi^2$-contrast is not symmetric and therefore no distance.
For arbitrary distributions it can be very large, i.e.
\[
	\chi^2(\nu,\mu) \leq \norm{\frac{\nu}{\mu}-1}{\infty},
	\quad \text{where} \quad
	\norm{\frac{\nu}{\mu}-1}{\infty}=\max_{x\in D} \abs{\frac{\nu(x)}{\mu(x)}-1}.
\]
From \cite[Theorem 3.3 p. 209]{bremaud} we have
\begin{equation} \label{conv_chi}
 \chi^2( \nu P^k,\pi ) \leq \beta^{2k} \; \chi^2 ( \nu, \pi),
\end{equation}
where $\beta=\max\set{\beta_1,\abs{\beta_{\abs{D}-1}}}$ denotes the second largest absolute value of the eigenvalues.
%Therefore we also have with \eqref{tv_chi} an upper bound for the total variation. 
Let us turn to another presentation of the convergence property.
We have
\begin{align*}
 \nu P^k(x)-\pi(x)	%=\sum_{y\in D} \nu(y) p^k(y,x) - \pi(x)
									&	= \sum_{y\in D} \frac{\nu(y)}{\pi(y)} p^k(y,x) \pi(y) - \pi(x) \\
									&	\underset{\text{rev.}}{=} \sum_{y\in D} \frac{\nu(y)}{\pi(y)} p^k(x,y) \pi(x) - \pi(x) \\
		&	= \sum_{y\in D} \frac{\nu(y)}{\pi(y)} p^k(x,y) \pi(x) - \sum_{y\in D} \frac{\nu(y)}{\pi(y)}\pi(y) \pi(x)\\
		&	= \sum_{y\in D} \frac{\nu(y)}{\pi(y)} (p^k(x,y)-\pi(y) )\pi(x).%=\sum_{y\in D} \frac{\nu(y)}{\pi(y)} (p^k(x,y)-\pi(y) )\pi(x).
\end{align*}
The second equality follows by the reversibility of the Markov chain.
For simplicity let%and a shorter notation
\[
d_k(x):= \sum_{y\in D} \frac{\nu(y)}{\pi(y)} (p^k(x,y)-\pi(y)),
\]
such that altogether 
\begin{align}
%\label{l_1}
%  \norm{l_k}{1} & =2\norm{\nu P^k(\cdot)-\pi(\cdot)}{\text{tv}} 
%  				\underset{\eqref{tv_chi}}{\leq} \min\set{2 , \beta^{k} \sqrt{\norm{\frac{\nu}{\pi}-1}{\infty}}},\\
\label{l_2}	
	\norm{d_k}{2} & = \sqrt{\chi^2( \nu P^k,\pi )} 
  				 \underset{\eqref{conv_chi}}{\leq} \beta^{k} \; \sqrt{\norm{\frac{\nu}{\pi}-1}{\infty}}.
\end{align}
 
Since $\beta<1$ we have an exponential decay of the norm with $k\to\infty$.
We define the weighted sequence spaces for $1\leq p \leq \infty$ by
\[
\ell_p=\ell_p(D,\pi):=\set{f \in \R^D: \norm{f}{p}^p=\sum_{x\in D} \abs{f(x)}^p \pi(x) < \infty}.
\]
It is clear that $\ell_p=\R^D$, since the state space has finite cardinality.
\begin{remark}
As we have seen the $\chi^2$-contrast corresponds to the $\ell_2$-norm of the function $d_k$.
Other tools for measuring the speed of convergence induce similar relations.
For instance
\[
\norm{d_k}{1}=2\norm{\nu P^k-\pi}{\text{tv}} \quad\text{and}\quad \norm{d_k}{\infty}=\norm{\frac{\nu P^k}{\pi}-1}{\infty}.
\]
The total variation corresponds to the $\ell_1$-norm of $d_k$ and the $\ell_\infty$-norm to the supremum-distance.%$\norm{\frac{\nu P^k}{\pi}-1}{\infty}$.
\end{remark}

\begin{remark}
The constant $\beta$ plays a crucial role in estimating the speed of convergence of the Markov chain to stationarity. In general it is not easy to handle $\beta_1$ or $\beta$, but there are different auxiliary tools, e.g. canonical path technique, conductance (see \cite{jerrum} and \cite{dia_stroock_eigen}), log-Sobolev inequalities and path coupling. For a small survey see \cite{randall}.
\end{remark}

\subsection{Norm of the transition matrix}
Let us consider $P$ and $S$ as operators acting on $\ell_p$. Then the functional $S$ maps arbitrary functions to constant functions. Let
\begin{align*}
\ell_p^0        := \ell_p^0( D,\pi)=\set{g\in \ell_p: S(g)=0} \quad \text{for}\quad 2\leq p\leq \infty. 
\end{align*}
The norm of $P$ as operator on $\ell_2^0$ and $\ell_4^0$ is essential in the analysis.
We state and show some results which are implied by %a classical theorem of interpolation theory, 
the Theorem of Riesz-Thorin. For a proof and an introduction we refer to \cite{bennett_sharpley}.

\begin{proposition}[Theorem of Riesz-Thorin] \label{riesz_thorin}
Let $1\leq p, q_1, q_2 \leq \infty$. Further let $\theta \in (0,1)$ and 
\[
\frac{1}{p}:=\frac{1-\theta}{q_1}+\frac{\theta}{q_2}
\]
and
\begin{align*}
& T:\ell_{q_1} \to \ell_{q_1}\quad \text{with} \quad \norm{T}{\ell_{q_1} \to \ell_{q_1}} \leq M_1,\\
& T:\ell_{q_2} \to \ell_{q_2}\quad \text{with} \quad \norm{T}{\ell_{q_2} \to \ell_{q_2}} \leq M_2.
\end{align*}
Then 
\[
\norm{T}{\ell_p\to \ell_p} \leq 2 M_1^{1-\theta} M_2^\theta.
\]
\end{proposition}  
Note that the factor two in the last inequality comes from the fact that we consider real valued functions $f$. 
In the following we show a relation between $P$, $P-S$ and $\beta$.
\begin{lemma} \label{P_S}
Let $P$ be a reversible transition matrix with respect to $\pi$ and $n\in\N$. Then
\begin{equation}
\label{P_S_2_ident}
   \norm{P^n-S}{\ell_2 \to \ell_2} = \norm{P^n}{\ell_2^0 \to \ell_2^0} = \beta^n.
\end{equation}
Furthermore if $\;2\leq p\leq \infty$ then
\begin{equation}
\label{P_S_p_ident}
   \norm{P^n}{\ell_p^0 \to \ell_p^0} \leq \norm{P^n-S}{\ell_p \to \ell_p}\leq 2. 
\end{equation}
\end{lemma}
\begin{proof}
The self-adjointness of $P$ implies $\norm{P}{\ell_2^0 \to \ell_2^0}=\max\set{\beta_1,\abs{\beta_{\abs{D}-1}}}=\beta$, such that $\norm{P^n}{\ell_2^0 \to \ell_2^0}  =  \beta^n$.
%At first we show \eqref{P_S_p_ident}. Simple transformations give
By
\begin{align*}
\norm{P^n-S}{\ell_2\to \ell_2} 	& =\sup_{\norm{f}{2}\leq1} \norm{(P^n-S)f}{2} = \sup_{\norm{f}{2}\leq1} \norm{P^n(f-S(f))}{2} \\
												& \leq \sup_{\norm{f}{2}\leq1} \sup_{\norm{g}{2}\leq1,\; S(g)=0} \norm{P^n g}{2} 
													= \norm{P^n}{\ell_2^0 \to \ell_2^0} 
\end{align*}	
and													
\begin{align*}
											\norm{P^n}{\ell_p^0 \to \ell_p^0} & = \sup_{\norm{g}{p}\leq1,\; S(g)=0} \norm{P^n g}{p}
												= \sup_{\norm{g}{p}\leq1,\; S(g)=0} \norm{P^n g-S(g)}{p}\\
												&	\leq \sup_{\norm{f}{p}\leq1} \norm{(P^n-S)f}{p} = \norm{P^n-S}{\ell_p\to \ell_p}
\end{align*}
claim \eqref{P_S_2_ident} and the first part of \eqref{P_S_p_ident} is shown. 
Finally, by applying the triangle inequality of the norm
\[
\norm{P^n-S}{\ell_p\to \ell_p} = \sup_{\norm{f}{p}\leq1} \norm{P^n f-Sf}{p} 
\leq \norm{P^n }{\ell_p\to \ell_p}+\norm{S}{\ell_p\to \ell_p}=2.
\]
\end{proof}
The next statement adds the result about the matrix norm which is used in the proof of the error bound.
\begin{lemma}
Let $P$ be a reversible transition matrix with respect to $\pi$ and $n\in\N$. Then 
\begin{equation}
	\norm{P^n}{\ell_4^0 \to \ell_4^0} \leq \, 2\sqrt{2}\;\beta^{n/2}.
	\label{norm_lp}
\end{equation} 
\end{lemma} 

\begin{proof}
By Lemma~\ref{P_S} we have 
\[
\norm{P^n-S}{\ell_2 \to \ell_2} = \beta^n \quad \text{and} \quad \norm{P^n-S}{\ell_\infty \to \ell_\infty} \leq 2.
\]
Then the result is an application of Proposition~\ref{riesz_thorin}, where $T=P^n-S$ and $q_1=2$, $q_2=\infty$, $p=4$ thus $\theta=\frac{1}{2}$. 
\end{proof}

%%%%%%%%%%%%%%%%%%%%%%%%%%%%%%%%%%%%%%%%%%%%%%%%%%%%%%%%%%%%%%%%%%%%%%%%
%%%%%%%%%%%%%%%%%%%%%%%%%%%%%%%%%%%%%%%%%%%%%%%%%%%%%%%%%%%%%%%%%%%%%%%%
%%%%%%%%%%%%%%%%%%%%%%%%%%%%%%%%%%%%%%%%%%%%%%%%%%%%%%%%%%%%%%%%%%%%%%%%

%\input{sections/error_bounds.tex}

\section{Error bounds} \label{err_sec}
In this section we mainly follow two steps to develop the error bound. At first a special case of method $S_{n,n_0}$ is considered. The initial distribution is the stationary one, thus it is not necessary 
to do a burn-in, i.e. $n_0=0$. % since the starting state is already chosen by the desired distribution. 
%to throw the first $n_0$ steps away since they are already chosen by the right distribution. 
Secondly we relate the result of the first step to the general case where the chain is initialized by a distribution $\nu$. The techniques which we will use are similar as in \cite{expl_error}.

\subsection{Starting from stationarity}
This is also called starting in equilibrium, i.e. the distribution of the Markov chain does not change, it is already balanced. In the following we will always denote $S_{n,0}$ as $S_n$. 
Let us start with stating and discussing a result from \cite[Prop. 2.1 p.3]{bassetti_dia},
which is similar to \cite[Theorem~1.9, p.~375]{lova_simo1}.
\begin{proposition} \label{expl_stat}
Let $f\in \R^D$.
Let $X_1,\dots,X_n$ be a reversible Markov chain with respect to $\pi$, given by $(P,\pi)$. Then

%chosen by $\pi$ and $X_1,\dots,X_n$ be a realization of the Markov chain which is given by a reversible $P$ as described above. Then 
\begin{equation} \label{err_present}
e_\pi(S_n,f)^2 = \frac{1}{n^2} \sum_{k=1}^{\abs{D}-1} 
\abs{a_k}^2 W(n,\beta_k),
						%-\frac{2}{n^2} \sum_{k=1}^{m-1} \abs{a_k}^2 \frac{\beta_k(1-\beta_k^{n})}{(1-\beta_k)^2}
%e(S_n,f)^2 \leq \left(\frac{c_1}{n} + \frac{c_2}{n^2}\right)\norm{f}{2}^2,  
\end{equation}
where
\begin{equation*}
 	a_k=\scalar{f}{u_k}_\pi\quad\mbox{and}\quad
 	W(n,\beta_k):=\frac{n(1-\beta_k^2)-2\beta_k(1-\beta_k^{n})}{(1-\beta_k)^2}.
\end{equation*}
\end{proposition}

\begin{proof}
 Let us consider $g:=f-S(f)\in\R^D$. Because of the orthogonal basis the presentation 
 $g(x)=\sum_{k=1}^{\abs{D}-1} a_k u_k(x)$ is given. The error obeys
 \begin{align*}
      e_\pi(S_n,f)^2 &	=\expect_{\pi,P}\abs{\frac{1}{n}\sum_{j=1}^n g(X_j )}^2
     						  =\frac{1}{n^2} \expect_{\pi,P} \abs{ \sum_{j=1}^n g(X_j)}^2 \\
     					&		=\frac{1}{n^2} \sum_{j=1}^n \expect_{\pi,P}\,g(X_j)^2
     								+ \frac{2}{n^2} \sum_{j=1}^{n-1} \sum_{i=j+1}^n \expect_{\pi,P}\,g(X_j)g(X_i).
 \end{align*}
 For $j\leq i$,
 \begin{align*}
  	 \expect_{\pi,P}\,g(X_i)g(X_j)	&= \sum_{k=1}^{\abs{D}-1} \sum_{l=1}^{\abs{D}-1} 
  																						a_ka_l\; \expect_{\pi,P}\, u_k(X_i)u_l(X_j)\\ 
 															 		%= \sum_{k=1}^{m-1} \sum_{l=1}^{m-1} a_ka_l\; \expect_{\pi,P}\, u_k(X_j)u_l(X_i)
 															&		\underset{\eqref{err_help}}{=} \sum_{k=1}^{\abs{D}-1} \sum_{l=1}^{\abs{D}-1} 
 																							a_ka_l\; \scalar{u_k}{P^{i-j}u_l}_\pi \\ 
 															&		= \sum_{k=1}^{\abs{D}-1} \sum_{l=1}^{\abs{D}-1} 
 																							a_ka_l\; \beta_l^{i-j}\scalar{u_k}{u_l}_\pi
 																	= \sum_{k=1}^{\abs{D}-1} a_k^2\; \beta_k^{i-j},
\end{align*} 		
 where the equality of the second line is due to the fact that the initial 
 step is chosen from the stationary distribution. 
 The last two equalities follow from the orthonormality of the basis of the eigenvectors.
 Altogether we have
 \begin{align*}
 		 e_\pi(S_n,f)^2 & =\frac{1}{n^2} \sum_{k=1}^{\abs{D}-1} a_k^2 
 				\left[ n + 2 \sum_{j=1}^{n-1} \sum_{i=j+1}^n \beta_k^{i-j} \right]\\
 		&	=\frac{1}{n^2} \sum_{k=1}^{\abs{D}-1} a_k^2 
 				\left[ n + 2 \frac{(n-1)\beta_k-n\beta_k^2+\beta^{n+1}_k}{(1-\beta_k)^2}\right]\\
 		&	=\frac{1}{n^2} \sum_{k=1}^{\abs{D}-1} 
\abs{a_k}^2 W(n,\beta_k)	.
  \end{align*} 

\end{proof}
%The error presentation above is rather complicated. 
Let us consider $W(n,\beta_k)$ to simplify and interpret Proposition~\ref{expl_stat}.
% presentation stated above in the stationary case.
\begin{lemma}
For all $n\in\N$ and $k\in\set{1,\dots,\abs{D}-1}$ we have
\begin{align}  \label{mono_inc}
W(n,\beta_k) & \leq W(n,\beta_1)\leq \frac{2n}{1-\beta_1}.
\end{align}
\end{lemma}

\begin{proof}
Let $x\in[-1,1)$, then we are going to show that $W(n,x)$ is monotone increasing, 
i.e. $W(n,\beta_k) \leq W(n,\beta_1)$. %in $x$ for $n\in\N$.
For $i\in\set{0,\dots,n-1}$ it is true that 
\[
x^{n-i}\leq1\quad \Longleftrightarrow \quad (1-x^i)\,x^{n-i} \leq 1-x^i \quad\Longleftrightarrow\quad x^{n-i}+x^i \leq 1+x^n.
\]
Therefore
\[
x^i+x^{i+1}+x^{n-i-1}+x^{n-i} \leq 2(1+x^n),
\]
and
\[
(1+x)\sum_{i=0}^{n-1} x^i =\frac{1}{2} \sum_{i=0}^{n-1} x^i+x^{i+1}+x^{n-i-1}+x^{n-i} \leq n(1+x^n).
\]
Now
\[
\frac{d \,W}{dx}(n,x)=-2\frac{(1+x)\sum_{i=0}^{n-1} x^i-n(1+x^n)}{(1-x)^2} \geq 0
\]
and the first inequality is shown. By
\[
W(n,x)\leq 	\begin{cases}
							\frac{n(1+x)-2xn}{1-x} 	& x\in[-1,0]\\
							\frac{n(1+x)}{1-x}			& x\in(0,1)
						\end{cases}
			\leq \frac{2n}{1-x}
\]
the claim is proven.
\end{proof}

An explicit formula of the error 
if the initial state is chosen by the stationary distribution is established. 
Let us discuss the worst case error of $S_n$.

\begin{proposition}  \label{prop_stat}
%Let $f\in \R^m$ and $S(f)=0$.
Let $X_1,\dots,X_n$ be a reversible Markov chain with respect to $\pi$, given by $(P,\pi)$. 
%\begin{enumerate}[(i)]
%\item \label{err_class}
Then
			\begin{equation}  \label{err_class_eq}
				\sup_{\norm{f}{2}\leq1}e_\pi(S_n,f)^2 = 
												\frac{1+\beta_1}{n(1-\beta_1)}-\frac{2\beta_1(1-\beta_1^n)}{n^2(1-\beta_1)^2}
												\leq \frac{2}{n(1-\beta_1)}.
			\end{equation}
%\item \label{spec_beta}
%			If $\beta_1\geq0$ then
%			\[
%			  \sup_{\norm{f}{2}\leq1}e(S_n,f)^2 \leq \frac{1+\beta_1}{n(1-\beta_1)}.
%			\]
%\end{enumerate}
\end{proposition} 

\begin{proof}
The individual error of $f$ is 
\begin{align*}
 e_\pi(S_n,f)^2	& 	\underset{\eqref{err_present}}{=} \frac{1}{n^2} \sum_{k=1}^{\abs{D}-1} \abs{a_k}^2 W(n,\beta_k)
								\leq 		\frac{\norm{f}{2}^2}{n^2} \max_{k=1,\dots, \abs{D}-1}  W(n,\beta_k) \\
	&%\label{stat_err}
			\underset{\eqref{mono_inc}}{=} \frac{\norm{f}{2}^2}{n^2} W(n,\beta_1)
								=\frac{1+\beta_1}{n(1-\beta_1)}\norm{f}{2}^2-\frac{2\beta_1(1-\beta_1^n)}{n^2(1-\beta_1)^2}\norm{f}{2}^2,
\end{align*}
where $a_k$ is chosen as in Proposition~\ref{expl_stat} and therefore $\sum_{k=1}^{\abs{D}-1} \abs{a_k}^2\leq\norm{f}{2}^2$.
%If $\beta_1\geq0$ the second term of the upper bound is damping the error such that 
%we can neglect it.
%in the worst case it vanishes completely. 
%After proving \eqref{spec_beta} let us turn to \eqref{err_class}. 
From the preceding analysis of the individual error we have an upper error bound. Now we consider $f=u_1$, where obviously $\norm{u_1}{2}=1$ and get by applying \eqref{err_present} that
\[
e_\pi(S_n,u_1)^2=\frac{1+\beta_1}{n(1-\beta_1)}-\frac{2\beta_1(1-\beta_1^n)}{n^2(1-\beta_1)^2}.
\] 
Thus the error bound is attained for $u_1$ and by \eqref{mono_inc} everything is shown.
\end{proof}

Finally an explicit presentation for the worst case error on the class of bounded functions with respect to $\norm{\cdot}{2}$ is shown. 
Notice, that \eqref{err_class_eq} is an equality, which means that the integration error is completely known if we start with the stationary distribution. %The assumption of \eqref{spec_beta} is always fulfilled if the underlying transition matrix is lazy (or positive semi-definite), since then all eigenvalues are positive. Here in contrast one just needs that the second largest eigenvalue is larger or equal to zero to get the upper error bound.
In some artificial cases this method even beats direct simulation, e.g. if one specific $\beta_i<0$ and the goal is to approximate $S(u_i)$ or if all eigenvalues are smaller than zero. 
In \cite[Remark 3, p.617]{eigen_neg} the authors state a simple transition matrix where $\beta_i=-\frac{1}{\abs{D}-1}$ for all $i$. Now one could think to construct a transition matrix where $\beta_1$ is close to $-1$ and therefore damp the integration error. But it is well known that this is not possible for large $\abs{D}$, since $\beta_1\geq -\frac{1}{\abs{D}-1}$.\\

In the next subsection we link the results to a more general framework, where the unrealistic assumption 
that the initial distribution is the stationary one is abandoned. 
  
\subsection{Starting from somewhere else}
In the next statement a relation between the error of starting by $\pi$ and the error
of starting not by the invariant distribution is established.

\begin{proposition} \label{connect_lem}
Let $f\in \R^D$ and $g:=f-S(f)$.
Let $X_1,\dots,X_{n+n_0}$ be a reversible Markov chain with respect to $\pi$, given by $(P,\nu)$. Then
%Let $P$ be a reversible transition matrix with stationary distribution $\pi$, let $X_1,X_2,\dots$ 
%be a Markov chain generated by $P$ with initial distribution $\nu$. 
%Then we get for $g:=f-S(f)\in \R^D$ that
\begin{align}  \label{connection}
%\expect_{\nu,P} \abs{S(f)-S_{n,n_0}(f)}^2 = & \expect_{\pi,P} \abs{S(f)-S_n(f)}^2\\
e_\nu(S_{n,n_0},f)^2 =  e_\pi(S_n,f)^2
+ \frac{1}{n^2}\sum_{j=1}^{n} L_{j+n_0}(g^2)
+ \frac{2}{n^2} \sum_{j=1}^{n-1} \sum_{k=j+1}^n L_{j+n_0}(gP^{k-j}g),
  %\int_ E \int_ E \frac{d\nu}{d\pi}(y) \;\left(K^{n_0+j}(x,dy)-\pi(dy)\right)\;g(x) P^{k-j}g(x)\;\pi(dx).
\end{align}
where
\[
L_i(h)= 	\sum_{x\in D} d_i(x) h(x) \pi(x) 
      = 	\sum_{x\in D} \sum_{y\in D} \frac{\nu(y)}{\pi(y)} (p^i(x,y)-\pi(y)) h(x) \pi(x).
\]
\end{proposition}

\begin{remark}
The proof of this identity is similar as in \cite{expl_error}, except for the fact that 
we study a finite state space and therefore integrals become sums. 
\end{remark}

\begin{proof}
It is easy to see, that
\begin{align*}
& \expect_{\nu,P} \abs{S(f)-S_{n,n_0}(f)}^2=\frac{1}{n^2} \sum_{j=1}^n \sum_{i=1}^n \expect_{\nu,P} (g(X_{n_0+j})g(X_{n_0+i}))\\
&= \frac{1}{n^2} \sum_{j=1}^n \sum_{x\in D} P^{n_0+j}g^2(x)\; \nu(x)+ \frac{2}{n^2} \sum_{j=1}^{n-1} \sum_{k=j+1}^n \sum_{x\in D} P^{n_0+j}(g P^{k-j}g)(x)\;\nu(x).
\end{align*}
For every function $h\in \R^D$ and $i\in\N$ under applying the reversibility the following transformation holds true
\begin{align*} 
& \sum_{x\in D} (P^i h)(x)\, \nu(x) 
	= \sum_{x\in D} \sum_{y\in D} h(y)\, p^i(x,y)\, \frac{\nu(x)}{\pi(x)}\, \pi(x)\\
& \underset{\text{rev.}}{=} 
	\sum_{x\in D} \sum_{y\in D} \frac{\nu(y)}{\pi(y)}  p^i(x,y)\,h(x) \, \pi(x)\\
& = \sum_{x\in D} h(x)\, \pi(x) 
	+ \sum_{x\in D} \sum_{y\in D} \frac{\nu(y)}{\pi(y)}  \left(p^i(x,y)-\pi(y)\right) h(x) \, \pi(x)\\
& \underset{\text{rev.}}{=} 
    \sum_{x\in D} (P^i h)(x)\pi(x) 
	+ \sum_{x\in D} \sum_{y\in D} \frac{\nu(y)}{\pi(y)}  \left(p^i(x,y)-\pi(y)\right) h(x) \, \pi(x).
\end{align*}
Using this in the setting above, formula \eqref{connection} is shown.
\end{proof}

Equation \eqref{connection} %and the result of Remark~\ref{connection_u_i} 
is still an error characterization where equality holds. 
We will estimate $L_k(h)$ to derive an upper bound.
%To derive an upper bound we need an estimation of $L_k(h)$. 
This depends very much on the 
speed of convergence from the chain to stationarity.% and also on the underlying measure 
%of this speed. 
%basically use total variation and $\chi^2$-contrast estimation of subsection~\ref{convergence}.% and total variation.

\begin{lemma}
Let  $h\in\R^D$, let again $\beta=\max\set{\beta_1,\abs{\beta_{\abs{D}-1}}}$. 
Then
\begin{align}
  \label{infty_1}
	\abs{L_k(h)} & %\leq 2 \norm{\frac{\nu}{\pi}-1}{\infty} 
		%\max_i \norm{\delta_i P^k (\cdot)-\pi(\cdot)}{\text{tv}} \cdot \norm{h}{1}
		\leq \beta^k  \sqrt{\norm{\frac{1}{\pi}}{\infty}}\sqrt{ \norm{\frac{\nu}{\pi}-1}{\infty} }  \cdot \norm{h}{1},
		\quad k\in\N,\\
  	\label{2_2}
	\abs{L_k(h)} &%\leq  \sqrt{\chi^2( \nu P^k,\pi )} \cdot \norm{h}{2}
	\leq \beta^k \sqrt{\norm{\frac{\nu}{\pi}-1}{\infty}}  \cdot \norm{h}{2},\quad k\in\N.
%	,\\[2ex]
%	 \label{1_infty}
%	\abs{L_k(h)} & \leq 2\norm{\nu P^k(\cdot)-\pi(\cdot)}{\text{tv}} \cdot  \norm{h}{\infty}
%	\leq \min\set{2,\beta^k \sqrt{\norm{\frac{\nu}{\pi}-1}{\infty}}}   \cdot \norm{h}{\infty}.
\end{align}

\end{lemma}
\begin{proof}
Let us consider $L_k(h)=\scalar{ d_k }{h}_\pi$.
After applying Cauchy-Schwarz inequality we obtain
\[
  \abs{L_k(h)} \leq \norm{d_k}{2}\norm{h}{2}.
\]
%for $1\leq p,q \leq \infty$ with $\frac{1}{p}+\frac{1}{q}=1$.
%If $p=1$ and $q=\infty$ then \eqref{1_infty} is immediately proven via \eqref{l_1}.
By applying \eqref{l_2} we showed \eqref{2_2}.
%For showing \eqref{2_2} we consider the case $p=q=2$ and get the claim with \eqref{l_2}.
Inequality \eqref{2_2} and $\norm{h}{2}\leq \sqrt{\norm{\frac{1}{\pi}}{\infty}} \norm{h}{1}$ imply \eqref{infty_1}.
\end{proof}
The ingredients for getting an explicit error bound for $S_{n,n_0}$ are gathered together. 
Mainly the last Lemma ensures an exponential decay of $L_k(h)$ 
which is used in the next Proposition. 
%Therefore, as we will see, by choosing the burn-in sufficiently large an optimal asymptotic behavior
%is reached.

\begin{proposition} \label{err_thm}
Let $X_1,\dots,X_{n+n_0}$ be a reversible Markov chain with respect to $\pi$, given by $(P,\nu)$.
Let $f\in \R^D$, $g:=f-S(f)$ and
\begin{align*}
V(\beta,n) & = \sum_{j=1}^n\beta^j + 2 \sum_{j=1}^{n-1} \sum_{k=j+1}^n \beta^k,\\
%=\frac{(2n-1)\beta^{n+2}-(2n+1)\beta^{n+1}+\beta^2+\beta}{(1-\beta)^2},\\
U(\beta,n) & = \sum_{j=1}^n\beta^j + 4\sqrt{2} \sum_{j=1}^{n-1} \sum_{k=j+1}^n \beta^{\frac{k+j}{2}}.
\end{align*}
\begin{enumerate}[(i)]
	\item \label{err_l_2}Then for $g\in \ell_2^0$ we have
	 \begin{align*}
	 e_\nu(S_{n,n_0},f)^2 \leq  e_\pi(S_n,f)^2
	 %\expect_{\pi,P} \abs{S(f)-S_n(f)}^2
	+\frac{V(\beta,n)}{n^2}\beta^{n_0} \sqrt{\norm{\frac{1}{\pi}}{\infty}}
	\sqrt{\norm{\frac{\nu}{\pi}-1}{\infty}} 
		  \norm{g}{2}^2.
	 \end{align*}

	\item \label{err_l_4}Then for $g\in \ell_4^0$ we have
	 \begin{align*}
	 e_\nu(S_{n,n_0},f)^2 \leq e_\pi(S_n,f)^2
	 %\expect_{\pi,P} \abs{S(f)-S_n(f)}^2
	 +\frac{U(\beta,n)}{n^2}\beta^{n_0} \sqrt{\norm{\frac{\nu}{\pi}-1}{\infty} } \norm{g}{4}^2.
	 \end{align*}

	\item \label{err_l_infty_improve}
	 Then for $g\in \ell_\infty^0$ we have %can improve the result of \eqref{err_l_2} such that
	 \begin{align*}
	 e_\nu(S_{n,n_0},f)^2 \leq e_\pi(S_n,f)^2
	 %\expect_{\pi,P} \abs{S(f)-S_n(f)}^2
	 +\frac{V(\beta,n)}{n^2}\beta^{n_0} \sqrt{\norm{\frac{\nu}{\pi}-1}{\infty} } \norm{g}{\infty}^2.
	%+\min\set{2,\frac{V(\beta,n)}{n^2}\beta^{n_0} \sqrt{\norm{\frac{\nu}{\pi}-1}{\infty} }} \norm{g}{\infty}^2.
	 \end{align*}
\end{enumerate}

\end{proposition}

\begin{proof}
	As we have seen in \eqref{connection} the error obeys
	\begin{align} \label{exact_error_a}  
	e_\nu(S_{n,n_0},f)^2 = e_\pi(S_n,f)^2
	% \expect_{\pi,P} \abs{S(f)-S_n(f)}^2\\
%& \notag\quad 
+ \frac{1}{n^2}\sum_{j=1}^{n} L_{j+n_0}(g^2)
+ \frac{2}{n^2} \sum_{j=1}^{n-1} \sum_{k=j+1}^n L_{j+n_0}(gP^{k-j}g).
\end{align}	
 Then by \eqref{infty_1}, Cauchy-Schwarz inequality and $\norm{P^{k-j}}{\ell_2^0\to \ell_2^0}=\beta^{k-j}$ we get
 \begin{align*}
 		\abs{L_{j+n_0}(g^2)} & 
 		\leq \sqrt{\norm{\frac{1}{\pi}}{\infty}} \sqrt{\norm{\frac{\nu}{\pi}-1}{\infty}}  \beta^{j+n_0} \norm{g}{2}^2 ,\\
 		\abs{L_{j+n_0}(gP^{k-j}g)}& 
 		\leq  \sqrt{\norm{\frac{1}{\pi}}{\infty}} \sqrt{\norm{\frac{\nu}{\pi}-1}{\infty}}  \beta^{k+n_0} \norm{g}{2}^2 .
 \end{align*}
 Putting this in the sums of equation \eqref{exact_error_a} and let $\e_0= \sqrt{\norm{\frac{1}{\pi}}{\infty}} \sqrt{\norm{\frac{\nu}{\pi}-1}{\infty}}  \beta^{n_0}$  we obtain
  \begin{align*}
 		\sum_{j=1}^n& \abs{L_{j+n_0}(g^2)} +
 		2\sum_{j=1}^{n-1} \sum_{k=j+1}^n \abs{L_{j+n_0}(gP^{k-j}g)}\\
 		&\leq \e_0 \norm{g}{2}^2  \sum_{j=1}^n\beta^j
 		+ \e_0 \norm{g}{2}^2 \sum_{j=1}^{n-1} \sum_{k=j+1}^n 2\beta^k\\
 		 &=  \e_0 \norm{g}{2}^2 \left( \sum_{j=1}^n\beta^j + \sum_{j=1}^{n-1} \sum_{k=j+1}^n 2\beta^k \right)
 		  = V(\beta,n) \cdot \e_0 \norm{g}{2}^2.
 \end{align*}
%The last equality follows from applying the finite geometric series several times.
%By multiplying of $1/n^2$ 
Thus claim \eqref{err_l_2} is shown.
Now we use \eqref{2_2} and 
\[
\norm{gP^{k-j}g}{2} \leq \norm{g}{\infty} \norm{P^{k-j}g}{2}  Ž
										\leq \norm{g}{\infty}^2 \norm{P^{k-j}}{\ell_2^0\to \ell_2^0} 
										\leq \norm{g}{\infty}^2 \beta^{k-j}
\]
to obtain
 \begin{align*}
 		\abs{L_{j+n_0}(g^2)} & 
 		\leq  \sqrt{\norm{\frac{\nu}{\pi}-1}{\infty}} \beta^{j+n_0}  \norm{g}{\infty}^2,\\
 		\abs{L_{j+n_0}(gP^{k-j}g)}& 
 		\leq  \sqrt{\norm{\frac{\nu}{\pi}-1}{\infty}} \beta^{k+n_0}  \norm{g}{\infty}^2 .
 \end{align*}
Exactly the same steps as in the proof of \eqref{err_l_2} follow, except for a different $\e_0=\sqrt{\norm{\frac{\nu}{\pi}-1}{\infty}} \beta^{n_0}$ and the supremum norm,
i.e. assertion \eqref{err_l_infty_improve} is proven.
Let us turn to \eqref{err_l_4}. Again we use \eqref{2_2} and estimate
\[
\norm{gP^{k-j}g}{2} \leq \norm{g}{4} \norm{P^{k-j}g}{4} \leq \norm{P^{k-j}}{\ell_4^0\to \ell_4^0} \norm{g}{4}^2 \underset{\eqref{norm_lp}}{\leq} 2\sqrt{2} \norm{g}{4}^2 \, \beta^{\frac{k-j}{2}}.
\]
Thus
 \begin{align*}
 		\abs{L_{j+n_0}(g^2)} & 
 		\leq  \sqrt{\norm{\frac{\nu}{\pi}-1}{\infty}} \beta^{j+n_0}  \norm{g}{4}^2,\\
 		\abs{L_{j+n_0}(gP^{k-j}g)}& 
 		\leq 2\sqrt{2}\; \sqrt{\norm{\frac{\nu}{\pi}-1}{\infty}} \beta^{\frac{k+j}{2}+n_0} \norm{g}{4}^2 .
 \end{align*}
For $\e_0=\sqrt{\norm{\frac{\nu}{\pi}-1}{\infty}} \beta^{n_0}$ we obtain
\begin{align*}
 		\sum_{j=1}^n& \abs{L_{j+n_0}(g^2)} +
 		4\sqrt{2}\sum_{j=1}^{n-1} \sum_{k=j+1}^n \abs{L_{j+n_0}(gP^{k-j}g)}\\
 		&\leq \e_0 \norm{g}{4}^2  \sum_{j=1}^n\beta^j 
 					+ \e_0 \norm{g}{4}^2 \sum_{j=1}^{n-1} \sum_{k=j+1}^n 4\sqrt{2}\, \beta^{\frac{k+j}{2}}\\
 		&=  \e_0 \norm{g}{4}^2 \left( \sum_{j=1}^n\beta^j + 4\sqrt{2} \sum_{j=1}^{n-1} \sum_{k=j+1}^n \beta^{\frac{k+j}{2}} \right)
 		= U(\beta,n) \cdot \e_0 \norm{g}{4}^2. \\
 \end{align*}
%The presentation of $U(\beta,n)$ follows by using the finite geometric series. 
Finally by substituting this in equation \eqref{exact_error_a} everything is shown.
\end{proof}
In the last Proposition we introduced $V(\beta,n)$ and $U(\beta,n)$.
These \-functions are bounded if $\beta<1$. 
By applying the infinite geometric series several times the following is proven.
\begin{lemma} \label{lemma_V}
For $n\in\N$ and $x\in[0,1)$ we have
\begin{alignat}{2}
%\label{V_n_U_n}		V(x,n) & \leq V(1,n) = n^2, 
%& \qquad U(x,n) &
%\leq U(1,n)=2\sqrt{2}n^2+(1-2\sqrt{2})n,\\
\label{V_x}		V(x,n) & \leq \frac{2}{(1-x)^2}, 
& \qquad U(x,n) & 
\leq \frac{4\sqrt{2}}{(1-x)(1-\sqrt{x})}.
\end{alignat}
%Especially 
%\begin{equation}
%\label{V_x_n} V(x,n)  \leq \frac{2n}{1-x}.
%\end{equation}
\end{lemma}
This implies that the asymptotic optimality is reached.

\subsection{Main Theorem} \label{main_thm_sec}
The following is the main result.
%After asking in the introduction for the presentation of the error bound, we answer this question by summarizing the attained results.

%After already developing an upper error bound of $S_{n,n_0}$ in Propostion~\ref{err_thm} we state the main theorem, which summarizes the attained results.

\begin{theorem} \label{main_thm}

Let $X_1,\dots,X_{n+n_0}$ be a reversible Markov chain with respect to $\pi$, given by $(P,\nu)$.
Let $f\in \R^D$ and $a_k=\scalar{f}{u_k}_\pi$.
Then
\begin{align*}
 \lim_{n\to \infty} n\cdot e_\nu(S_{n,n_0},f)^2 &=\lim_{n\to \infty} n\cdot e_\pi(S_{n},f)^2 = \sum_{k=1}^{\abs{D}-1} \abs{a_k}^2 \frac{1+\beta_k}{1-\beta_k}.
\end{align*}
\begin{enumerate}[(i)]
	\item \label{expl_err_l2} If we consider $f\in \ell_2$ then
	 \begin{align*}
    e_\nu(S_{n,n_0},f)^2
		&\leq \frac{2}{n(1-\beta_1)}\norm{f}{2}^2 
				%-\frac{2\beta_1(1-\beta_1^n)}{n^2(1-\beta_1)^2}
  			+				 
				 \frac{2\sqrt{\norm{\frac{1}{\pi}}{\infty}} \sqrt{\norm{\frac{\nu}{\pi}-1}{\infty}} 
				 \beta^{n_0} }{n^2(1-\beta)^2}\norm{f}{2}^2.
	 \end{align*}
	\item \label{expl_err_l4} If we consider $f\in \ell_4$ then
	 \begin{align*}
e_\nu(S_{n,n_0},f)^2
		&\leq \frac{2}{n(1-\beta_1)}\norm{f}{4}^2 
				%-\frac{2\beta_1(1-\beta_1^n)}{n^2(1-\beta_1)^2}
  			+ \frac{16\sqrt{2}\sqrt{\norm{\frac{\nu}{\pi}-1}{\infty}}\beta^{n_0}}{n^2(1-\beta)(1-\sqrt{\beta})} \norm{f}{4}^2.
	 \end{align*}
	\item \label{expl_err_infty}
	 If we consider $f\in \ell_\infty$ then
	 \begin{align*}
e_\nu(S_{n,n_0},f)^2
		&\leq \frac{2}{n(1-\beta_1)}\norm{f}{\infty}^2
		+ \frac{4%\sqrt{\norm{\frac{1}{\pi}}{\infty}} 
		\sqrt{\norm{\frac{\nu}{\pi}-1}{\infty}} 
				 \beta^{n_0} }{n^2(1-\beta)^2}\norm{f}{\infty} ^2	
				  \end{align*}
\end{enumerate}

\end{theorem}
\begin{proof}
By \eqref{connection} and the fact that the remaining terms are going quadratic to zero as $n$ goes to infinity, we see that the asymptotic result holds true.
%it is immediately seen, that a certain burn-in does not effect the asymptotic optimality. 
%It is obvious 
For $f\in \ell_2$ we have $\norm{f-S(f)}{2}\leq \norm{f}{2}$ and furthermore
if $p\not=2$ then
\[
\norm{f-S(f)}{p}\leq \norm{f}{p}+\abs{S(f)} \leq \norm{f}{p}+\norm{f}{1} \leq 2\norm{f}{p}.
\]
Thus, via Proposition~\ref{err_thm}, Proposition~\ref{prop_stat} and Lemma~\ref{lemma_V} everything is shown.
\end{proof}

Notice, that from the estimate of Proposition~\ref{err_thm} it follows immediately that
\[
\lim_{n\to\infty} n\cdot e_\nu(S_{n,n_0},f)^2 \leq \lim_{n\to\infty} n\cdot e_\pi(S_{n,n_0},f)^2 \leq \frac{1+\beta_1}{1-\beta_1} \norm{f}{2}^2. 
%\leq \frac{1+\beta_1}{1-\beta_1}
\] 
Thus there is no gap between the estimate and the asymptotical behavior. 
Also notice, that the upper bounds are continuous in the sense that if the initial distribution $\nu$ is $\pi$ then we obtain the bound of Proposition~\ref{prop_stat}. The dependence of the bounds of \eqref{expl_err_l4} and \eqref{expl_err_infty} in Theorem~\ref{main_thm} on the initial distribution is encouraging for an extension to general state spaces.
(For MCMC on general state spaces we refer to \cite{rosenthal,novak, expl_error}.)
%For examples of Markov chains and estimates of $\beta_1$ or $\beta$ we also refer to \cite{novak,small_world}.)
But the dependence of the initial distribution on the estimate in the $\ell_2$-case is disillusioning because of the additional factor of $\norm{\frac{1}{\pi}}{\infty}$.\\   
 
In \cite[Theorem 8, p.10]{expl_error} a similar $\ell_\infty$-bound of $S_{n,n_0}$ for general state spaces is developed. This result holds for lazy, reversible Markov chains and may also be applied in the present setting, i.e. if the state space is finite. 
In \cite{expl_error} the asymptotic error limit is not attained. 
Thus we could improve the error bound and weaken the laziness condition, i.e. it is enough that $\beta_1=\beta$. In \cite[Thm. 12.19, p.165]{peres} the authors obtained for another error term a comparable bound where the chain starts deterministically. Very recently in \cite{niemiro} a similar result concerning the integration error for $f\in \ell_\infty$ was shown where the Markov chain is not necessarily reversible.

%%%%%%%%%%%%%%%%%%%%%%%%%%%%%%%%%%%%%%%%%%%%%%%%%%%%%%%%%%%%%%%%%%%%%%%%
%%%%%%%%%%%%%%%%%%%%%%%%%%%%%%%%%%%%%%%%%%%%%%%%%%%%%%%%%%%%%%%%%%%%%%%%
%%%%%%%%%%%%%%%%%%%%%%%%%%%%%%%%%%%%%%%%%%%%%%%%%%%%%%%%%%%%%%%%%%%%%%%%

%\input{sections/burn_in.tex}

\section{Burn-in} \label{burn_in_sec}

Let us assume that computer resources for the MCMC method for $N$ time steps are available, i.e. $N=n+n_0$.
We want to choose the burn-in $n_0$ and the number of $n$ such that the error bound is as small as possible. The burn-in $n_0$ should be large but this implies that $n$ is possibly quite small depending on how much resources we have. On the other hand $n$ should be large which again implies that $n_0$ is possibly small. There is obviously a trade-off between choosing the parameters. In the next statement we consider the error for an explicitly given burn-in, where for simplicity $\beta_1=\beta$.

\begin{corollary}
\label{main_coro}
Let $f\in\R^{D}$ be given and let 
\[
n_0= \max\set{\left\lceil \frac{\log(C)}{\log(\beta^{-1})}\right\rceil,0}.
\] 
\begin{enumerate}[(i)]
\item Let $C=\sqrt{\norm{\frac{1}{\pi}}{\infty}} \sqrt{\norm{\frac{\nu}{\pi}-1}{\infty}} $, then 
\[
\sup_{\norm{f}{2}\leq1}e_\nu(S_{n,n_0},f)^2
		\leq \frac{2}{n(1-\beta)} 
				+				 
				 \frac{2}{n^2(1-\beta)^2}.
\]
\item Let $C=16\sqrt{2}\sqrt{\norm{\frac{\nu}{\pi}-1}{\infty}}$, then 
\[
\sup_{\norm{f}{4}\leq1}e_\nu(S_{n,n_0},f)^2
		\leq \frac{2}{n(1-\beta)} 
				+				 
				 \frac{1}{n^2(1-\beta)(1-\sqrt{\beta})}.
\]
\item Let $C=2\sqrt{\norm{\frac{\nu}{\pi}-1}{\infty}}$, then 
\[
\sup_{\norm{f}{\infty}\leq1}e_\nu(S_{n,n_0},f)^2
		\leq \frac{2}{n(1-\beta)} 
				+				 
				 \frac{2}{n^2(1-\beta)^2}.\]
\end{enumerate}
\end{corollary}

Note, that in the $\ell_\infty$- and $\ell_2$-case the error bound is the same. Just the constant $C$ which comes in by the density is different. This suggestion of the burn-in is justified in the following.

\subsection{Numerical experiments}
Suppose $C$ (very large), $\beta$ (close to one) and resources $N$ are given.
The worst case error 
for $\norm{f}{2}\leq1$ or $\norm{f}{\infty}\leq1$ is bounded by
\[
b_\infty(n,n_0):=\sqrt{\frac{2}{n(1-\beta)}+\frac{2C \beta^{n_0}}{n(1-\beta)^2}}
\]
and if we consider $\norm{f}{4}\leq1$ it is bounded by
\[
b_4(n,n_0):=\sqrt{\frac{2}{n(1-\beta)}+\frac{C \beta^{n_0}}{n(1-\beta)(1-\sqrt{\beta})}}.
\]
Since $N=n+n_0$ we can compute with a numerical procedure (here using Maple) the optimal choice of the burn-in denoted by $n^4_{\text{opt}}$, $n^\infty_{\text{opt}}$ to minimize the upper error bounds. (This is a simple one dimensional minimization problem with different parameters.)

\begin{table}[htb]
\begin{tabular}{|c|c|c|c|c|}

\hline 

&&&& \\[-2ex]

$N$	& $\beta$		&  $n^4_{\text{opt}}$	& $n^\infty_{\text{opt}}$	& $n_0=\left\lceil \log(C)/\log(\beta^{-1}) \right\rceil$ 

\\ 			& 			& (by Maple)	& (by Maple)			& (suggested above)    		\\[1ex]  		
\hline & & & & \\[-1.5ex] 

$10^4$ & $0.9$   	&   $656$				& $656$ 					& $656$		\\
$10^5$ & $0.9$		&	  $656$				&	$656$						&	$656$ 	\\
$10^4$ & $0.99$   &   $6867$			& $6867$ 					&	$6873$	\\
$10^5$ & $0.99$		&	  $6873$			&	$6873$					&	$6873$ 	\\
$10^4$ & $0.999$  &   $8001$			& $8001$ 					&	$69043$  \\
$10^5$ & $0.999$	&	  $68977$			&	$68977$				&	$69043$ 	\\[1ex]
\hline
\end{tabular} 
\\[1.5ex]
\caption{For $C=10^{30}$ where $n^i_{\text{opt}}$ 
minimizes $b_i(N-n^i_{\text{opt}},n^i_{\text{opt}})$,\; $i=4,\infty$.}
 \label{tab_burn_in}
\end{table}

Table~\ref{tab_burn_in} gives a collection of typical results. It turned out that the above suggested lower bound is close to the optimal choice. The computed value $n^4_{\text{opt}}$ and $n^\infty_{\text{opt}}$ is almost the same as $n_0=\left\lceil \log(C)/\log(\beta^{-1}) \right\rceil$. In the case $N=10^4$ and $\beta=0.999$ Theorem~\ref{main_thm} gives for no choice of $n$ and $n_0$ an error smaller than one. 

\begin{figure}[htb] 
  \begin{center}
    \includegraphics[height=9cm]{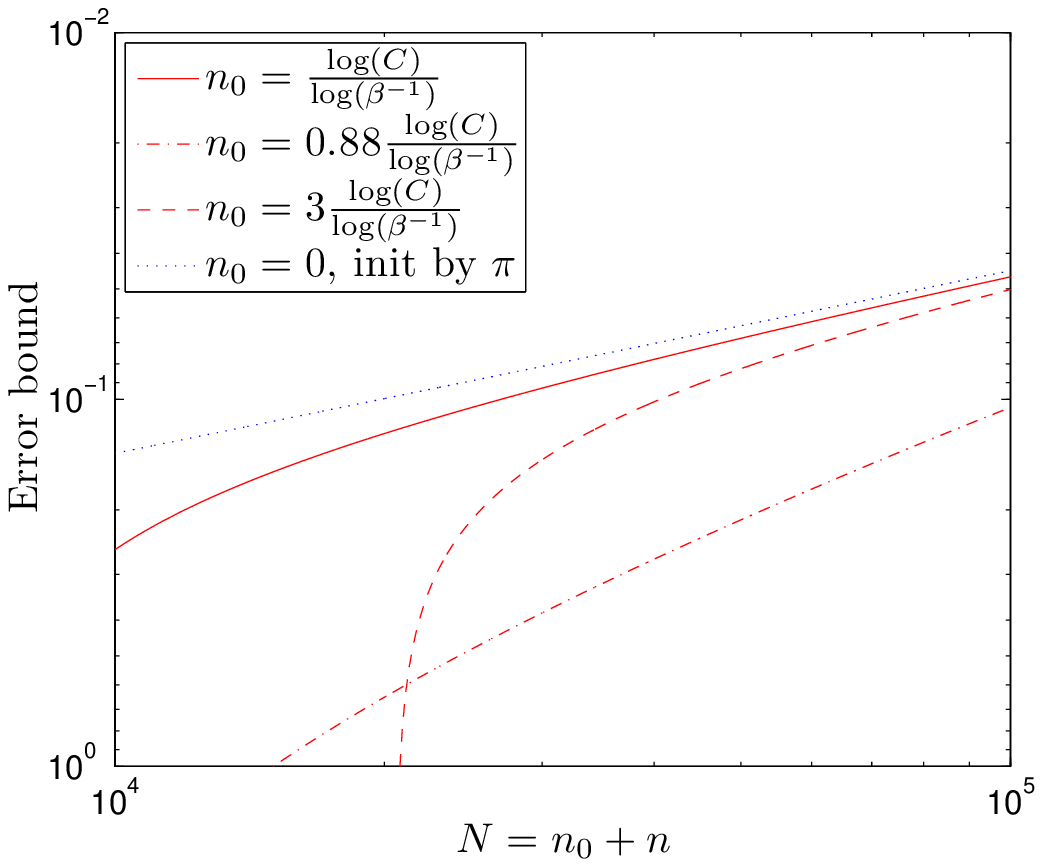}
    \caption{For $\beta=0.99$ and $C=10^{30}$.}
    \label{diff_asymp}
    \end{center}  
\end{figure}

For different $n_0$ we plotted in Figure~\ref{diff_asymp} 
\[
b_4(N-n_0,n_0) \quad \text{and} \quad 
e_\pi(S_N,u_1)=\sqrt{\frac{1+\beta_1}{N(1-\beta_1)} -\frac{2\beta_1(1-\beta_1^N)}{N^2(1-\beta_1)^2}}.
\]
%$$  and $e_\pi(S_N,u_1)=\sqrt{\frac{1+\beta_1}{N(1-\beta_1)} -\frac{2\beta_1(1-\beta_1^N)}{N^2(1-\beta_1)^2}}$.
Roughly spoken one may see in Figure~\ref{diff_asymp} that if the burn-in is chosen too small a vertical shifting takes place and if the burn-in is chosen to large a horizontal shifting takes place. 
The asymptotic behavior is the same, i.e. for the long run the error of $S_{n,n_0}$ converges to the error of $S_n$.
If $\beta$ and $C$ are given we chose the burn-in as suggested above. 
If there is an estimate of $\log(C)/\log(\beta^{-1})$ one should ensure that it is not smaller than the real ratio. 
As seen in Figure~\ref{diff_asymp} if it is slightly smaller there is already strong influence. By choosing the burn-in too large the influence is less heavy. 
 
\begin{figure}[htb] 
  \begin{center}
    \includegraphics[height=9cm]{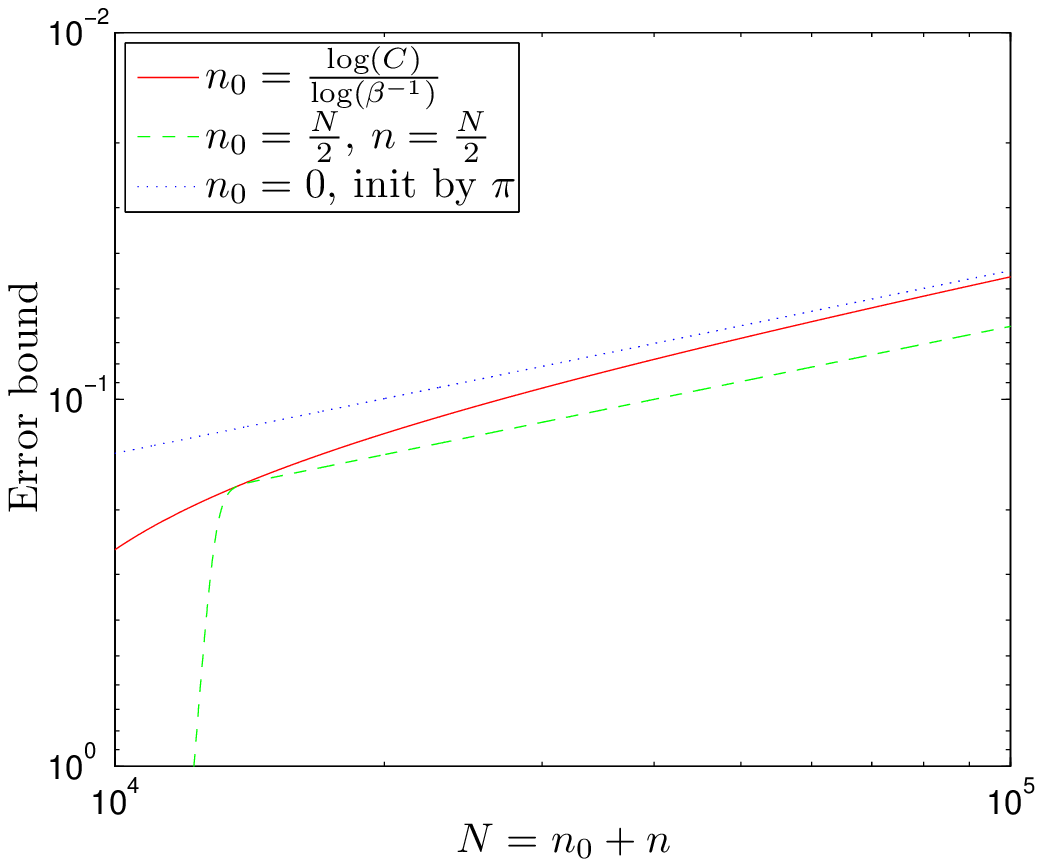}
    \caption{For $\beta=0.99$ and $C=10^{30}$.}
    \label{exact}
    \end{center}  
\end{figure}

Finally if there is no estimation or computation of the parameters $\beta$ or $C$ a 
simple but very efficient strategy is given by choosing $n=n_0=\frac{N}{2}$ (for even $N$). 
In Figure~\ref{exact} we see $b_4(\frac{N}{2},\frac{N}{2})$, $b_4(N-n_0,n_0)$ and $e_\pi(S_N,u_1)$.
%This has mainly the advantage that we do not need any information about $\beta$ or $D$. 
In the asymptotic behavior we pay the price of a factor of $\sqrt{2}$, i.e. the asymptotic error is $\sqrt{2}$ times larger than $e_\pi(S_N,u_1)$ where we started in equilibrium. 
This strategy works well and reaches the same convergence rate as choosing the burn-in as suggested above,
which is seen in Figure~\ref{exact}.

\section*{Acknowledgements}
The author thanks Erich Novak and Aicke Hinrichs for their valuable comments.  

%%%%%%%%%%%%%%%%%%%%%%%%%%%%%%%%%%%%%%%%%%%%%%%%%%%%%%%%%%%%%%%%%%%%

%*******************************%

%*****************************************************%
% Please choose one of the following options:         %
%                                                     %
% VERSION A is for BIBTeX application.                %
% You have to add the names of your databases:        %
%                                                     %
%%\bibliography{database}
%                                                     %
% IF you prefer VERSION B for the LaTeX standard      %
% bibliography environment, please use the same style %
% as produced by mcma.bst (see the examples in the    %
% instructions for authors).                          %
%                                                     %
%\begin{thebibliography}{99}                          %
%  \bibitem{Lange:2007}                               %
%  T.~Lange and I.~E. Shaparlinski, Distribution of   %
%  some sequences of points on elliptic curves,       %
%  \emph{J. Math. Cryptol.} \textbf{1} (2007), 1--11. %
%\end{thebibliography}                                %
%                                                     %
%*****************************************************%

\end{document}